\documentclass[12pt]{amsart}

\usepackage{amscd}
\usepackage{amsfonts}
\usepackage{amssymb}
\usepackage{amsmath}

\newcommand{\la}{\lambda}
\newcommand{\al}{\alpha}
\newcommand{\be}{\beta}

\newcommand{\f}{\varphi}
\newcommand{\h}{\phi}

\newcommand{\abs}[1]{\vert #1\vert}
\newcommand{\ov}{\overline}
\newcommand{\Ll}{\mathcal{L}}

\newcommand{\GL}{\mathop{\rm GL}\nolimits}

\newcommand{\ra}{{\rightarrow}}
\newcommand{\lra}{{\longrightarrow}}
\newcommand{\Null}{\mathop{\rm Null}\nolimits}
\newcommand{\Aut}{\mathop{\rm Aut}\nolimits}

\newcommand{\ISO}{\mathop{\rm ISO}\nolimits}

\newcommand{\CC}{\mathbb{C}} 
 
\newcommand{\RR}{\mathbb{R}}
\newcommand{\ZZ}{\mathbb{Z}}
\newcommand{\NN}{\mathbb{N}}

\numberwithin{equation}{section}

\newtheorem*{ta}{Theorem A}
\newtheorem*{cb}{Corollary B.}
\newtheorem*{tc}{Theorem C}
\newtheorem*{cd}{Corollary D}
\newtheorem{pr}{Proposition}[section]
\newtheorem{co}{Corollary}[section]
\newtheorem{lm}{Lemma}[section]

\theoremstyle{definition}
\newtheorem{de}{Definition}[section] 
\newtheorem{re}{Remark}[section]

\begin{document}

\title[Transformations of compact l.c.K. manifolds]
{Transformations of compact locally conformally K\"ahler manifolds}
\author{Yoshinobu Kamishima}
\address{Department of Mathematics, Tokyo Metropolitan 
University,\newline
Minami-Ohsawa 1-1, Hachioji, Tokyo 192-0397, Japan}
\email{kami@comp.metro-u.ac.jp}

\author{Liviu Ornea}
\address{University of Bucharest, Faculty of Mathematics\newline
14 Academiei str.,
70109 Bucharest, Romania}
\email{lornea@imar.ro}
\date{\today}
\begin{abstract}We characterize compact locally conformally
K\"ahler (l.c.K.) manifolds
by means of the existence of a purely conformal, holomorphic circle action.
As an application, we determine the structure of the
compact locally conformally
K\"ahler  manifolds with parallel Lee form.
We introduce  the Lee-Cauchy-Riemann (LCR) transformations as a class of
diffeomorphisms preserving the specific $G$-structure of l.c.K.
manifolds. Then we characterize the
Hopf manifolds, up to holomorphic isometry, as
compact l.c.K. manifolds admitting a certain closed LCR action of 
$\CC^*$.
\end{abstract}

\maketitle

\section{Introduction}

Let $(M,g,J)$ be a connected, complex Hermitian manifold of complex 
dimension $n\geq 2$. We denote its 
 fundamental two-form $\omega$ defined by $\omega(X,Y) =g(X,JY)$.
\begin{de}
If $\omega$ satisfies
the integrability condition
\begin{equation*}
d\omega=\theta\wedge\omega\quad \text{with}\; d\theta=0
\end{equation*}
the manifold is  called \emph{locally conformally K\"ahler} (l.c.K.).
\end{de}
The closed one-form
$\theta$ is called the Lee form and it encodes the geometric properties 
of
such a manifold (see \cite{DO} and the bibliography therein). Note that 
in
complex dimension at least $3$, if $d\omega=\theta\wedge\omega$,
then $\theta$ is automatically closed. The second condition in the 
definition
is necessary only on complex surfaces.
Let $\nabla^g$ be the Levi-Civita connection of the l.c.K. metric $g$.
The Weyl connection
$$D=\nabla^g-\frac{1}{2}\{\theta\otimes Id+Id\otimes\theta-g\otimes 
\theta^\sharp\}$$
is almost complex ($DJ=0$) and preserves the conformal class of $g$
($Dg=\theta\otimes g$). Hence, a locally conformally K\"ahler manifold 
is
a Hermitian-Weyl manifold. The converse is also true in complex 
dimension at least $6$ (see \cite{PPS}).

A locally conformally K\"ahler manifold  $(M,g,J)$ whose
Lee form is $\nabla^g$-parallel is called a \emph{Vaisman manifold}. 
This subclass was introduced by I. Vaisman under
the name of generalized Hopf
manifolds; but not all of the Hopf surfaces have parallel Lee form,
cf. \cite{Be}, \cite{GO}. The structure of compact Vaisman manifolds is
 better understood. Here we recall some facts.
At the topological level, it is known that $b_1$ is odd (\cite{Va1}).
Moreover, the Lee field $\theta^\sharp$ is real analytic
($\Ll_{\theta^\sharp}J=0$) and
$g$-Killing. The leaves of the foliation
$\mathop{\Null}\, \theta$ have an induced Sasakian structure.
We note that  $\theta^\sharp$ and $J\theta^\sharp$ generate a
$2$-dimensional, holomorphic, integrable distribution whose orthogonal
distribution is preserved by $J$.
When $\theta^\sharp$ generates an
$S^1$ action by $g$-isometries,
its action is quasi-regular ({i.e.,} without fixed points),
and so the quotient space is a Sasakian orbifold whose
characteristic field is the projection of $J\theta^\sharp$.
The Riemannian universal covering space
$\tilde M$ of $M$ then splits as $\RR\times \tilde N$
where $\tilde N$ is a Sasakian manifold
(as concerns Sasakian manifolds, see \cite{Bl}).
An immediate example (see \cite{Va}) 
is the Hopf manifold $\CC^n-\{0\}/\Gamma$,
where $\Gamma$ is the cyclic group generated by $(z^i)\mapsto (2z^i)$,
endowed with the projection of the metric
$(\sum{\abs{z^i}}^2)^{-1}\sum dz^i\otimes d\ov{z}^i$,
globally conformal with
the standard flat metric of $\CC^n$.
 As the Hopf manifold is diffeomorphic
with $S^1\times S^{2n-1}$ it cannot bear any K\"ahler metric.
Moreover, the Lee
form of this metric is easily seen to be parallel with respect to
the Levi-Civita connection.

In this paper we  try to understand the influence of
the structure of some transformations groups  on the geometry of a
compact l.c.K. manifold. 

We first consider $Aut_{l.c.K.}(M)$, the (compact) group of all conformal,
holomorphic diffeomorphisms. Some of its elements may have the
stronger property  to not restrict to a local isometry on any chart of
the l.c.K. structure; we call them \emph{purely conformal}. If 
 a circle in $Aut_{l.c.K.}(M)$ is purely conformal, 
its lift to the universal
cover is an $\RR$ action (cf. Lemma 2.1). As a consequence, 
a compact, semi-simple subgroup of $Aut_{l.c.K.}(M)$ 
cannot contain purely conformal transformations (Corollary \eqref{croc}).
Our first main 
result (Theorem A) states that
the existence of a purely conformal circle in $Aut_{l.c.K.}(M)$
assures the existence of a metric with parallel Lee form in the given
conformal class. Note that, till now,  l.c.K. manifolds with parallel Lee 
form has been characterized
in the l.c.K. class by means of second order curvature conditions
(involving the Ricci curvature) or, 
when the Lee
field is regular, as circle bundles over Sasakian manifolds (as
described above). From this result we derive the structure theorem for 
Vaisman manifolds (without any assumption on the regularity of its
natural foliations). Roughly speaking, a Vaisman manifold is isometric
with  $\displaystyle S^1\mathop{{\times}}_{Q}^{} W$ where 
$\displaystyle S^1\ra M\stackrel {\pi}{\lra} W/Q$ is a
Seifert fiber space over a Sasakian orbifold $W/Q$ (see Corollary B).  

In the second half of the paper we introduce a larger group of
diffeomorphisms, containing the l.c.K.-ones, the
\emph{Lee-Cauchy-Riemann} (LCR) transformations. 
These are characterized as preserving the specific $G$-structure of a
l.c.K. manifold. Precisely, with respect to an orthonormal  coframe 
$\{\theta, \theta\circ J, \theta^\al,\ov{\theta}^\al\}_{\al=1,
\cdots, n-1}$ adapted to a l.c.K. manifold $(M,g,J)$, a LCR
transformation $f$ acts as:
 \begin{equation*}
\begin{split}
f^*\theta=\theta,\ \  f^*(\theta\circ J)=\lambda\cdot(\theta\circ 
J),& \\
f^*\theta^{\al}=\sqrt \lambda\cdot\theta^\be U^{\al}_{\be}+
(\theta\circ J)\cdot v^{\al},& \\
f^*{\bar \theta}^{\al}=\sqrt \lambda\cdot
{\bar\theta}^\be\ov{U}^\al_\be+(\theta\circ J)\cdot\ov{v}^{\al}, &
\end{split}
\end{equation*}
$\lambda$ being a positive, smooth function  and
$U^\al_\be$ a matrix in ${\rm U}(n-1)$. The main results of this part
(Theorem C and Corollary D), essentially say that a compact 
l.c.K. manifold with the
 Lee field generated by a circle action and which admits 
a closed, non-compact flow of 
LCR-transformations is 
bi-holomorphically isometric  with a finite quotient of a Hopf manifold. 
The interest (and difficulty) of this statement lies 
in the consideration of a non-compact subgroup of the group of LCR
transformations. The proof relies heavily on techniques specific to
CR geometry. It is to be noted that, to the authors' knowledge, the
only previous attempts to characterize the Hopf manifolds among the
Vaisman manifolds used curvature (precisely: conformally flatness) or
spectral properties (see 
\cite{DO}).

\section{Locally conformally K\"ahler transformations}

Let $(M,J)$ be a complex manifold.
A l.c.K. structure on $(M,J)$ is a local family
$\{U_\al, g_\al\}_{\al\in\Lambda}$ where 
$\{U_\al\}_{\al\in \Lambda}$ is an open cover and
$g_\al$ is a K\"ahler metric on each $U_\al$.
On non-empty intersections $U_\al\cap U_\be$,
there exist positive
constants $\la_{\be\al}$
such that $g_\be=\la_{\be\al}g_\al$.
It turns out that $\{\la_{\be\al}\}$ is a $1$-cocycle
on $M$.
Viewed $H^1(M;\RR^+)$ as
the sheaf cohomology of locally defined smooth positive functions,
there exists a
local family $\{f_\al,U_\al\}_{\al\in\Lambda}$ consisting of
smooth positive functions defined on each neighborhood
such that 
$\displaystyle\delta^0f(\al,\be)=\frac{f_\al}{f_\be}=\lambda_{\be\alpha}$ on
$U_\al\cap U_\be\neq\emptyset$.
Setting $g|U_\al=f_\al\cdot g_\al$, there exists a Hermitian metric 
$g$
on $M$ which is locally conformal to K\"ahler metrics.
\begin{de}
Two l.c.K. structures
$\{U_\al, g_\al\}_{\al\in\Lambda}$,
$\{{U}_\al, {g'}_\al\}_{\al\in\Lambda}$
are equivalent if there exists
a local family $\{c_\al, U_\al\}_{\al\in\Lambda}$
(more precisely, a refinement)
where each $c_\al $ is a constant number such that
${g'}_{\al}=c_\al\cdot g_{\al}$.
\end{de}

Given a l.c.K. manifold $(M,g,J)$,
since the Lee form $\theta$ is exact locally, there exists
an open cover $\{U_\al\}_{\al\in \Lambda}$ and positive functions
$f_\al$ defined
on $U_\al$ such that $df_\al=\theta|{U_\al}$. It follows that
the local metrics defined on $U_\al$ by 
$g_\al=e^{-f_\al}\cdot g{|U_\al}$ are K\"ahler and conformally related:
$e^{f_\be}\cdot g_\be=e^{f_\al}\cdot g_\al$ on  non-empty intersections
$U_\al\cap U_\be$. Consequently, $\la_{\be\al}=e^{f_\al- f_\be}$
are positive
constants satisfying the cocycle condition.
Hence  there is  a locally conformally K\"ahler structure
$\{U_\al, g_\al\}_{\al\in\Lambda}$ adapted to
each l.c.K. manifold $(M,g,J)$.

\begin{pr} \label{lck}
The set of equivalence classes of
 l.c.K. structures on a complex manifold $(M,J)$
is in one-to-one correpondence
with the set of conformal classes of l.c.K. metrics on $(M,J)$.
\end{pr}
\begin{proof}
Supose that a l.c.K. structure
$\{U_\al, g_\al\}_{\al\in\Lambda}$
is equivalent to  another $\{{U}_\al, g'_\al\}_{\al\in\Lambda}$.
By the definition, $g'_{\al}=c_\al\cdot g_{\al}$.
As $g'_{\be}={\lambda'}_{\be\al}g'_\al$,
${g}_{\be}={\lambda}_{\be\al}{g}_\al$, we have
$\lambda'_{\be\al}=c_\be\cdot\lambda_{\be\al}\cdot c_{\al}^{-1}$
on the intersection $U_\al\cap U_\be$.
Let $\displaystyle\delta^0f(\al,\be)=\frac{f_\al}{f_\be}=\lambda_{\be\alpha}$ on
$U_\al\cap U_\be\neq\emptyset$. Similarly for
$\displaystyle\lambda'_{\be\al}=\frac{f'_\al}{f'_\be}$ as above.
Thus there is a global function
$\tau$ on $M$ such that
$\tau|U_\al=c_\al f'_\al f_\al^{-1}$.
By the definition, the l.c.K. metrics $g$, $g'$ satisy
$g|U_\al=f_\al\cdot g_\al,\ g'|U_\al=f'_\al\cdot g'_\al$.
Hence $\tau\cdot g|U_{\al}=g'|U_\al$ for each $\al$.
Thus the equivalence class of $\{U_\al, g_\al\}_{\al\in\Lambda}$
determines the conformal class  $[g]$.

Conversely,
if $g'=\lambda\cdot g$ for
a l.c.K. metric $g$ where $\lambda$
is some positive function, note that
the fundamental two-form $\omega'$ satisfies 
$d\omega'=\theta'\wedge \omega'$ where $\theta'=\theta+d{\rm log}\lambda$.
As $d\omega'=0$,  $(M,g',J)$ is also a l.c.K. manifold.
Let $df_\al=\theta{|U_\al}$ and
$df'_\al={\theta'}{|U_\al}$ as above.
Since $df'_\al={\theta'}{|U_\al}=\theta{|U_\al}+d{\rm log}\lambda$,
there is some constant $c_\al$ such that
${\rm log}c_\al+f'_\al=f_{\al}+{\rm log}\lambda$
on $U_\al$. In particular,
$\displaystyle e^{f'_\al}e^{-f_\al}=\lambda c_\al^{-1}$.
By the defintion,
$\displaystyle g'_\al=e^{-f'_\al}\cdot g'=e^{-f'_\al}\cdot c_\al\cdot g=c_\al\cdot g_\al$.
So the induced l.c.K. structures
$\{U_\al, g_\al\}_{\al\in\Lambda}$,
$\{{U}_\al, {g'}_\al\}_{\al\in\Lambda}$
are equivalent.
\end{proof}

Let us denote $Aut_{l.c.K.}(M)$
the group of l.c.K. transformations: 
it consists of all diffeomorphisms of $M$
preserving the l.c.K. structure adapted to $g$.
Explicitly, $\f\in Aut_{l.c.K.}(M)$
satisfies the conditions:
\begin{gather*}
\f_*J=J\f_*\\
\f^*g_\be=\mu_{\be\al}g_\al, \quad \mu_{\be\al}=ct.>0 \quad 
\text{whenever}\;
\f(U_\al)\subseteq U_\be
\end{gather*}
As regards $Aut_{l.c.K.}(M)$, the following result was proved
by the first named author in \cite{Ka}.
\begin{pr}
Let $(M,g,J)$ be a locally conformally K\"ahler manifold of complex 
dimension
$n\geq 2$. Then  $Aut_{l.c.K.}(M)$ is a closed subgroup of the group of 
all
conformal diffeomorphisms of $(M,g)$. Moreover, if $M$ is
compact, $Aut_{l.c.K.}(M)$ is a compact Lie group.
\end{pr}

By the definition, the first statement is clear. As for the second one, 
if $Aut_{l.c.K.}(M)$
is noncompact, then also the group of all conformal transformations of 
$(M,g)$ is noncompact.
Then, by the celebrated result of Obata and Lelong-Ferrand, $(M,g)$ is 
conformally equivalent
with the sphere $S^{2n}$, $n\geq 2$. Hence $M$ is simply connected, thus 
K\"ahlerian.
It is well known that a conformally flat K\"ahler manifold of dimension
 $n\geq 3$ is flat, which is impossible, while
$S^{4}$ has no complex structure.

Since  a l.c.K. manifold $(M,g,J)$ is compact in our case,
$Aut_{l.c.K.}(M)$ is a compact Lie group.
Averaging the metric $g$ by the compact group, we obtain
a $Aut_{l.c.K.}(M)$-invariant metric $g'$ conformal to $g$.
\begin{pr}
Given a  l.c.K. structure on a compact complex manifold $(M,J)$,
there exists a l.c.K. metric $(M,g,J)$
for which the group $Aut_{l.c.K.}(M)$ acts as isometries.
Especially, $Aut_{l.c.K.}(M)$ leaves invariant the Lee form $\theta$ and
 the anti-Lee form
$\theta\circ J$.
\end{pr}

 We are interested in studying the possible subgroups of
$Aut_{l.c.K.}(M)$. It is necessary to distinguish a special
class of
conformal transformations:
\begin{de}
An element $a\in Aut_{l.c.K.}(M)$ is called \emph{purely conformal} if
$a^*g_\al=c_{\al\be}g_\be$ with $c_{\al\be}\neq 1$ for all $\al, \be$
(\emph{i.e.} $a$ is not 
a local isometry with respect to any of  the local K\"ahler metrics).  
A subgroup of $Aut_{l.c.K.}(M)$ is called \emph{purely conformal}
if it contains only
purely conformal  holomorphic transformations.
\end{de}

\begin{re}\label{ex}
On any compact locally conformally K\"ahler manifold with parallel Lee 
form, if
the Lee field generates a circle action, this is
purely conformal.
We shall see this phenomenon more rigorously in Corollary B.
( See also \cite{DO}, \cite {Va1}.)
Such examples of l.c.K. manifolds with
regular (respectively  quasi-regular) parallel Lee field are
 principal circle bundles over
 compact Hodge manifolds (respectively  Brieskorn manifolds of a certain type).
\end{re}
In order to see the action of the fundamental group,
we are looking at the universal covering space $\tilde M$
rather than to $M$ itself in this paper.
Let  $(\tilde M, \tilde g)$ be the universal Riemannian covering space
of $(M,g)$  and denote equally with $J$ the lifted complex structure. 
Then 
the lifted Lee form $\tilde \theta$ is  exact (as $\tilde M$ is simply 
connected): 
$\tilde \theta=d\tau$. Thus
$(\tilde M, \tilde g, J)$ is globally conformally K\"ahler. Let
$h$ be the K\"ahler metric $h=e^{-\tau}\cdot \tilde g$. 
Locally, on the inverse 
image $\tilde U_\al$ of $U_\al$, one has
\begin{equation}\label{po}
h=e^{-\tau}\cdot\tilde f_\al\tilde g_\al
\end{equation}
Note that, as $h$ and $\tilde g_\alpha$ are K\"ahler metrics,
the function  $e^{-\tau}\cdot\tilde f_\al$ is a constant.
Clearly $\pi_1(M)$ acts by holomorphic $\tilde g$-isometries, hence by
holomorphic conformal transformations with respect to $h$.
But it is well-known  that a holomorphic conformal transformation
of a K\"ahler manifold is homothetic. The converse costruction is also
possible. We obtain:

\begin{pr}$($\cite{Va1}$)$ Let $(M,J)$ be a complex manifold and $\tilde 
M$ its
universal Riemannian cover. Then $M$ admits a locally conformally 
K\"ahler
structure if and only if $\tilde M$ admits
a K\"ahlerian structure with respect to which $\pi_1(M)$ acts by 
homothetic
holomorphic transformations.
\end{pr}

The next lemma provides a characterization of purely conformal l.c.K. 
transformations:
\begin{lm} \label{ess}
A locally conformally K\"ahler circle action on $M$ lifts to an
$\RR$-action by
holomorphic, non-trivial 
 $h$-ho\-mo\-the\-ties on $\tilde M$ if and only if it is purely 
conformal.
\end{lm}
\begin{proof}
A l.c.K. $S^1$ action lifts to an $S^1$ or a $\RR$-action by homotheties 
on $\tilde M$. Let  us show that, owing to the purely conformal 
character of $S^1$ in $Aut_{l.c.K.}(M)$, it is impossible to lift to an 
$S^1$-action.
 Let $b\in S^1$ be such that, for a given $\al$, $b(U_\al)\subseteq 
U_\al$ (such a $b$ must be small enough,  close to $1$.) Then
$b^*g_\al=c_\al\cdot g_\al$ with $c_\al=const.\neq 1$ on $U_\al$.
This $b$ lifts to a $\tilde b$
acting on $\tilde M$ and we have
$\tilde b^*\tilde g_\al=c_\al\cdot \tilde g_\al$ on $\tilde U_\al$.  
The
K\"ahler metric $h$ is related to the lifted globally conformally 
K\"ahler
metric $\tilde g$ by $h=e^{-\tau}\tilde g$, hence \eqref{po} implies 
(recall that
$e^{-\tau}\cdot \tilde f_\al$ is constant on $U_\al$):
$$\tilde b^*h= (e^{-\tau}\cdot \tilde f_\al)\cdot \tilde b^*\tilde 
g_\al=
(e^{-\tau}\cdot \tilde f_\al)\cdot c_\al\tilde g_\al.$$
Now, if $\tilde b$ is an isometry, then $c_\al=1$, contradiction.
Thus, if $\rho$ is the homomorphism
 defined on the set of all $h$-homotheties of
$\tilde M$ with values in $\RR^+$ by the formula $t^*h =\rho(t)h$, we 
derive
$\rho(\tilde b)\neq 1$. Hence, the lift of $S^1$ cannot be $S^1$, 
otherwise
its image by $\rho$ would be a compact \emph{non-trivial} subgroup of 
$\RR^+$, contradiction.
\end{proof}

If $g'=\lambda\cdot g$, then the K\"ahler metric is obtained
as  $\displaystyle h'=e^{\tau'}\cdot \tilde {g}'$ where
 $\displaystyle\tilde \theta'=d\tau'$.
As $h,h'$ are K\"ahler,
we have $h'=c'\cdot h$ where
$c'=e^{-\tau'}\tilde \lambda e^{\tau}$ is a constant.
So the lifted group $\RR$ consists of nontrivial $h$-homotheties
if and only if it consists of nontrivial $h'$-homotheties.
\begin{co}
A purely conformal circle action does not depend on the conformal class
of l.c.K. metrics $[ g ]$.
\end{co}
The following result imposes a restriction 
on the existence of purely conformal circle actions
of a l.c.K. manifold:
\begin{co} \label{croc}
A  compact semisimple subgroup of $Aut_{l.c.K.}(M)$ 
cannot contain any purely conformal transformation.
\end{co}
\begin{proof}
Suppose, \emph{ad absurdum},  that $M$ admits a purely conformal $S^1$ 
action. According
to Lemma \ref{ess} (and to its proof),
this action lifts to an $\RR$-action by $h$-homotheties and
$\rho$ is injective. But $\rho$ can be understood both as being defined 
on $\RR$ and on   $\pi_1(M)$, because
$\RR$ and $\pi_1(M)$ are both contained in the group of
$h$-homotheties of $\tilde M$.
 Further, for a fixed $x\in M$, consider the evaluation map
$ev:G\rightarrow M$, $ev(a)=ax$. We denote also by $ev$ the 
restriction of
this map
to the considered $S^1$, so obtain the commutative diagram:
\begin{center}
\begin{picture}(100,60)(-10,-40)
\put(-25,0){$G$}
\put(-12,4){\vector(1,0){70}}
\put(20,10){$ev$}
\put(60,0){$M$}
\put(-23,-45){$S^1$}
\put(-8,-39){\vector(2,1){70}} 
\put(25,-28){$ev$}  
\put(-20,-35){\vector(0,1){30}} 
\end{picture}
\end{center}
where the vertical arrow is injective. Taking into account the
previous observation, at the homotopy level we then
have the commutative diagram:
\begin{center}
\begin{picture}(110,60)(-10,-40)
\put(-25,0){$\pi_1(G)$}
\put(10,4){\vector(1,0){40}}
\put(25,10){$ev_\sharp$}
\put(60,0){$\pi_1(M)$}
\put(100,4){\vector(1,0){40}}
\put(120,10){$\rho$}
\put(150,0){$\RR^+$}
\put(-23,-45){$\ZZ$}
\put(-20,-35){\vector(0,1){30}}
\put(-10,-42){\vector(4,1){150}}
\put(-15,-39){\vector(2,1){70}}
\put(25,-25){$ev_\sharp$}
\put(60,-35){$\rho$}
\end{picture}
\end{center}

Hence, $\rho\circ ev_\sharp (\ZZ)$ is 
infinite in $\RR^+$.
But, $G$ being semi-simple, it
has no torus as a direct summand
and has finite $\pi_1(G)$. Thus, chasing on the upper side of the diagram, 
the image of
$\ZZ$ in $\RR^+$ through 
$\rho\circ ev_\sharp$
is finite. This contradiction completes the proof.
\end{proof}

 Our first main result shows that the existence of purely conformal
circles in $Aut_{l.c.K.}(M)$
characterizes the existence of metrics with parallel Lee form in a
given conformal class:
\begin{ta}
Let $(M,g,J)$ be a compact locally conformally K\"ahler, non-K\"ahler
manifold. If $Aut_{l.c.K.}(M)$ contains 
a purely conformal subgroup $S^1$, then there exists a metric with 
parallel Lee form in the
conformal class of $g$.
\end{ta}
\subsection*{Proof of Theorem A}\hfill

\par Let us call $\xi$ the generating vector field of the
$\RR$-action on $\tilde M$
 and let $<\f_t>$ be the one-parameter subgroup generated by
$\xi$.
Using the fact that $\f_t$ are holomorphic maps we derive
\begin{equation}\label{doi}
[\xi, J\xi]=\lim_{t\rightarrow 0}\frac{\f_{t *}(J\xi)-J\xi}{t}=
\lim_{t\rightarrow 0}\frac{J\f_{t *}(\xi)-J\xi}{t}=0,
\end{equation}
thus the one-parameter subgroups generated by $\xi$ and $J\xi$ commute.

\begin{lm} \label{compl}
The one-parameter subgroup $<\psi_t>$ generated by $J\xi$ is global.
\end{lm}
\begin{proof}
As the deck group of the covering $\pi$ preserves $\xi$
and commutes with $J$, $J\xi$ projects to a one-parameter subgroup
on $M$. This one is global, by compactness of $M$ and lifts to a global
one-parameter subgroup on $\tilde M$ which clearly coincides 
with $\psi_t$.
\end{proof}
\begin{re}
As $\xi$ is a real analytic vector field, so is $J\xi$ (cf. \cite{Ko},
p. 76).
\end{re}
\par Let $\Omega$ be the K\"ahler form of the K\"ahler metric $h$ on
$\tilde M$.  (Note that $\Omega=e^{-\tau}\cdot \tilde \omega$.)
Because $\RR$ acts  by
$h$-homotheties, we can write:
 \begin{equation}\label{rho}
 \f_t^*\Omega =\rho(t)\cdot \Omega, \quad t\in \RR, \,
\rho(t)\in \RR^+.
\end{equation}
As $\rho$ is a nontrivial, continuous homomorphism, $\rho(t)=e^{at}$
for some constant $a\neq 1$. We may normalize $a=1$ so that
\eqref{rho} is described as
 \begin{equation}\label{rho1}
 \f_t^*\Omega =e^t\cdot \Omega, \quad t\in \RR.
\end{equation}
\begin{lm} \label{s1}
Let $s:\tilde M\rightarrow \RR$ be the smooth map defined as
$s(x)=\Omega(J\xi_x,\xi_x)$. Then $1$ is a regular value of $s$, hence
$s^{-1}(1)$ is a codimension $1$, smooth submanifold of $\tilde M$.
\end{lm}
\begin{proof}
Note that
\begin{equation*}
\begin{split}
s(\f_tx)&=\Omega(J\xi_{\f_tx},\xi_{\f_tx})=
\Omega(\f_{t*}J\xi_x,\f_{t*}\xi_x)\; \text{by \eqref{doi}}\\
&=e^t\Omega(J\xi_x,\xi_x)=e^ts(x)\; \text{by \eqref{rho1}}
\end{split}
\end{equation*}
As $s(x)\neq 0$, $s^{-1}(1)\neq \emptyset$. On the other hand, from
$$\Ll_\xi\Omega=
\lim_{t\rightarrow 0}\frac{\f_t^*\Omega-\Omega}{t}=\Omega,$$
we obtain for $x\in s^{-1}(1)$:
$$ds(\xi_x)=\xi s(x)=(\Ll_\xi s)(x)=(\Ll_\xi\Omega(J\xi,\xi))(x)=
\Omega(J\xi_x,\xi_x)=s(x)=1.
$$
This proves that $ds:T_x\tilde M\rightarrow \RR$ is onto and
$s^{-1}(1)$ is a codimension-$1$, smooth submanifold of $\tilde M$.

\end{proof}

Let now $W=s^{-1}(1)$. We can prove:
\begin{lm}
$W$ is connected. The evaluation map $H:\RR\times W\rightarrow \tilde 
M$,
defined by $H(t, w)=\f_tw$ is an equivariant diffeomorphism.
\end{lm}
\begin{proof}
Let $W_0$ be a component of $s^{-1}(1)$ and let $\RR\cdot W_0$ be the 
set
$\{\f_tw\; ;\; w\in W_0, t\in \RR\}$. As $\RR$ acts freely and
$s(\f_tx)=e^ts(x)$, we have $\f_tW_0\cap W_0=\emptyset$ for $t\neq 
0$.
Thus $\RR\cdot W_0$ is an open subset of $\tilde M$. We now prove that
it
is also closed. Let $\ov{\RR\cdot W_0}$ be the closure of $\RR\cdot
W_0$ in $\tilde M$. We choose a limit point $p=\lim \f_{t_i}w_i\in
\ov{\RR\cdot W_0}$. Then $s(p)=\lim s(\f_{t_i}w_i)=\lim
e^{t_i}s(w_i)=\lim e^{t_i}$. We denote $t=\log s(p)$, $t=\lim 
t_i$,
so $\f_t^{-1}(p)=\lim \f_{t_i}^{-1}(\lim \f_{t_i}w_i)=\lim w_i$.
Since $s^{-1}(1)$ is regular (\emph{i.e.} closed with respect to the
relative topology induced from $\tilde M$), so is its component $W_0$.
Hence $\f_t^{-1}p\in W_0$. Therefore $p=\f_t(\f_t^{-1}p)\in \RR\cdot
W_0$, proving that $\RR\cdot W_0$ is closed in $\tilde M$. In
conclusion, $\RR\cdot W_0=\tilde M$. Now, if $W_1$ is another
component of $s^{-1}(1)$, the same argument shows
 $\RR\cdot W_1=\tilde M$. Thus, $\RR\cdot
W_0=\RR\cdot W_1$. This implies $W_0=W_1$, in other words $W$ is
connected.
\end{proof}

\begin{lm} \label{re}
$W$ has a contact form $\eta$ for which $\pi_*J\xi$ is the Reeb
$($characteristic$)$ field.
\end{lm}

\begin{proof}
Let $i:W\rightarrow \tilde M$ be the inclusion and $\pi:\tilde
M\rightarrow W$ be the canonical projection. For $X_w\in T_wW$ we put
\begin{equation}\label{trei}
\eta_w(X_w)=e^{-t}\Omega(\tilde X_{\f_tw},\xi_{\f_tw})
\end{equation}
where $\tilde X_{\f_tw}\in T_{\f_tw}\tilde M$ such that $\pi_*\tilde
X=X$. As $\RR\rightarrow \tilde M\stackrel{\pi}{\rightarrow}W$ is a
fiber bundle with $T\RR=<\xi>$, $\eta$ is well-defined. To prove that
$\eta$ is a contact form, we first note that by the definition
\begin{equation}\label{et}
\pi^*\eta=-e^{-t}\cdot\iota_\xi\Omega \;\, \text{on each $T_{\f_tw}\tilde 
M$}.
\end{equation}
Then
$$i^*\pi^*\eta=(\pi i)^*\eta
=\eta=i^*(-e^{-t}\cdot\iota_\xi\Omega)=-i^*\iota_\xi\Omega\;\, 
\text{on $W$}.
$$
Now recall that $\Omega=\Ll_\xi\Omega=d\iota_\xi\Omega+\iota_\xi
d\Omega=d\iota_\xi\Omega$, thus
$$d\eta=-di^*\iota_\xi\Omega=-i^*d\iota_\xi\Omega=-i^*\Omega$$
and, moreover,
$$\eta(\pi_*(J\xi_w))=\Omega(J\xi_w,\xi_w)=s(w)=1\;\, 
\text{by the definition}.$$
Hence, $\eta\wedge d\eta^{n-1}\neq 0$ on $W$ showing that $\eta$ is a
contact form.

Let us now show that $\pi_*J\xi$ is the characteristic field of $\eta$. For any
distribution $D$ on $\tilde M$, denote $D^\perp$ the orthogonal 
distribution
with respect to the
metric $h$. Then $\pi_*:\xi^\perp\rightarrow TW$ is an isomorphism and
induces an isomorphism
$\pi_*:\{\xi,J\xi\}^\perp\rightarrow\mathop{\Null}\, \eta$.
So, $\Omega(\tilde X,\xi)=h(\tilde
X,J\xi)=0$, $\Omega(\tilde X,J\xi)=h(\tilde
X,-\xi)=0$ for $\tilde X\in  \{\xi,J\xi\}^\perp$. We now show that
$d\eta(\pi_*J\xi, X)=0$ for any $X\in \mathop{\Null}\, \eta$. We have
$2d\eta(\pi_*J\xi, X)=-\eta([\pi_*J\xi,X])$. Let $\pi_*\tilde X=X$ 
for some $\tilde X\in \{\xi, J\xi\}^\perp$. Then, using \eqref{et}
on $W$ (\emph{i.e.} for $t=0$):
$$2d\eta(\pi_*J\xi, X)=-\pi^*\eta([J\xi,\tilde X])=
\iota_\xi\Omega([J\xi,\tilde X])=\Omega(\xi,[J\xi,\tilde 
X]).$$
But $\displaystyle[\tilde X,\xi]=\mathop{\lim}_{t\rightarrow 0}\frac 
1t(\tilde
X-\f_{-t*}\tilde X)$ and
$$\Omega(\f_{-t*}\tilde X,J\xi)=h(\f_{-t*}\tilde 
X,\xi)=e^{-t}h(\tilde
X,\f_{t*}\xi)=e^{-t}h(\tilde X,\xi)=0$$
because $\tilde X\in \{\xi,J\xi\}^\perp$. So:
$$\Omega([\tilde X,\xi],J\xi)=\lim_{t\rightarrow 0}\frac{\Omega(\tilde
X, J\xi)-\Omega(\f_{-t*}\tilde X, J\xi)}{t}=0.$$
Since
\begin{equation*}
\begin{split}
3d\Omega(\tilde X, \xi,J\xi)&=\tilde X\Omega(\xi,J\xi)-
\xi\Omega(\tilde X,J\xi)+(J\xi)\Omega(\tilde X,\xi)\\
&-\Omega([\tilde X,\xi],J\xi)
-\Omega([\xi,J\xi],\tilde X)-\Omega([J\xi,\tilde X],\xi)=0,
\end{split}
\end{equation*}
we have
\[
\tilde X\Omega(\xi,J\xi)-\Omega([J\xi,\tilde X],\xi)=0.
\]
If $\exp:T_w\tilde M\ra \tilde M$ is the exponential map with respect to
$h$, then the subset
$\{t,\ \exp(\xi^{\perp})\}$ constitutes a coordinate neighborhood of 
$\tilde M$ at $w$. In particular, the above $\tilde X$
is a linear combination of coordinate vector fields around $w$
without containing $\frac {\partial}{\partial t}$.
On the other hand, noting
$\Omega(J\xi_{\varphi_t w},\xi_{\varphi_t w})=s({\varphi_t w})
=e^{t}$ on $\tilde M$ from \eqref{s1}, it follows
$\tilde X\Omega(\xi,J\xi)=-\tilde X\cdot e^t=0$.
We finally deduce
$\Omega([J\xi,\tilde X],\xi)=0$ on $W$ and the proof of the lemma is
complete.

\end{proof}
It is now easy to verify the following:
\begin{co}
$(W,i^*h,\pi_*J\xi)$ is a Sasakian manifold. In particular, if we denote 
by $\phi$ the tangent part of $J$  to $W$ (hence
$\phi=\nabla^{h'}\pi_*J\xi$), we can say that
$(\eta, \phi)$ is a pseudo-Hermitian structure on $W$.
\end{co}

We have derived from \eqref{et} that
\begin{equation}\label{eta}
d(e^{t}\pi^*\eta)=-\Omega \;\, \text{on}\  \tilde M.
\end{equation}For the given l.c.K. metric $g$, 
the K\"ahler metric $h$ is obtained as $h=e^{-\tau}\cdot \tilde g$ where
$d\tau=\tilde \theta$. (See \eqref{po}.)
As $\omega$ is the fundamental two-form of $g$, note
that $\Omega=e^{-\tau}\cdot \tilde\omega$.

Put
\begin{equation}\label{bar}
\bar\Theta=-e^{-t}\cdot d(e^{t}\pi^*\eta)\ (=+e^{-t}\cdot\Omega).
\end{equation}Then
the Hermitian metric $\bar g(X,Y)=\bar\Theta(JX,Y)$
is a l.c.K. metric of $\RR\times W$
on which $\RR$ acts as isometries.
It is obvious that $\bar g$ has
the parallel Lee form $-dt$.
We obtain that
\begin{equation}\label{co}
\bar\Theta=\mu\cdot \tilde\omega\ \text{(equivalently}\
\bar g=\mu\cdot\tilde g)
\end{equation}where  $\mu=e^{-(t+\tau)}:
\tilde M\ra \RR^+$ is a smooth map.

\begin{lm}\label{inv}
$\pi_1(M)$ acts by isometries of $\bar g$.
\end{lm}
\begin{proof}
We prove the following two facts:
\begin{enumerate}
\item $\gamma^*\pi^*\eta=\pi^*\eta$ for every $\gamma\in\pi_1(M)$.
\item  $\gamma^*e^t=\rho(\gamma)\cdot e^t$ where $\rho:\pi_1(M)\ra \RR^+$
is a homomorphism similar to the case of $\RR=\{\f_\theta\}$.
\end{enumerate}
First note that as $\RR=\{\f_\theta\}$ centralizes $\pi_1(M)$,
$\gamma_*\xi=\xi$ for $\gamma\in\pi_1(M)$.
Recall (cf. \eqref{re}) that
$\pi_*:\xi^\perp\rightarrow TW$,
$\pi_*:\{\xi,J\xi\}^\perp\rightarrow\mathop{\Null}\, \eta$
are isomorphic. As $\pi_1(M)$
 acts on $\tilde M$ as  holomorphic homothetic transformations (cf. 
Proposition 2.4), $\pi_1(M)$ leaves $X\in\{\xi,J\xi\}^\perp$ invariant.
If $X\in\{\xi,J\xi\}^\perp$, then $\gamma^*\pi^*\eta(X)
=\eta(\pi_*\gamma_*X)=0$. As $\pi_*J\xi$ is a characteristic vector 
field for $\eta$, $\gamma^*\pi^*\eta(J\xi)
=\eta(\pi_*\gamma_*J\xi)=\eta(\pi_*J\xi)=1$. This concludes that
$\gamma^*\pi^*\eta=\pi^*\eta$. On the other hand,
if we note $\gamma_*\xi=\xi$,
then
\begin{eqnarray*}
\gamma^*(\iota_\xi\Omega)(X)&=&\Omega(\xi,\gamma_*X)=
\Omega(\gamma_*\xi,\gamma_*X)\\
&=&\gamma^*\Omega(\xi,X)=\rho(\gamma)\cdot\Omega(\xi,X)\\
&=&\rho(\gamma)\cdot \iota_\xi\Omega(X)
\end{eqnarray*}where $\rho(\gamma)$ is a positive constant number.
As $\gamma^*\pi^*\eta=\pi^*\eta$ from {\bf 1.} and
$\pi^*\eta=-e^{-t}\cdot\iota_\xi\Omega$ from \eqref{et},
we obtain that $\gamma^*e^{-t}\cdot\rho(\gamma)=e^{-t}$.
Equivalently, $\gamma^*e^{t}=\rho(\gamma)\cdot e^{t}$.
This shows {\bf 1} and {\bf 2}.
Now, we prove \eqref{inv}. 

From \eqref{bar},
\begin{eqnarray*}
\gamma^*\bar\Theta&=&\gamma^*(-e^{-t}\cdot d(e^{t}\pi^*\eta))\\
 &=&-\rho(\gamma)^{-1}\cdot e^{-t}d(\rho(\gamma)\cdot
e^{t}\gamma^*\pi^*\eta)\\
&=&-e^{-t}\cdot d(e^{t}\pi^*\eta)=\bar\Theta.
\end{eqnarray*}
Since $\bar g(X,Y)=\bar\Theta(JX,Y)$ by the definition,
$\pi_1(M)$ acts through holomorphic isometries of $\bar g$.

\end{proof}

From this lemma,  the covering map $p:\tilde M\ra M$ induces a
l.c.K. metric $\hat g$ with parallel Lee form $\hat\theta$ 
on $M$ such that $p^*\hat g=\bar g$ and $p^*\hat\theta=-dt$.
Since $\tilde g$ is a lift of $g$ to $\tilde M$,
using the equation  \eqref{co}, we derive
\[
\gamma^*\mu\cdot\gamma^*\tilde g=\gamma^*\mu\cdot\tilde g
=\gamma^*\bar g=\bar g=\mu\cdot\tilde g
\] therefore 
$\gamma^*\mu=\mu$.
Since $\mu$ factors through a map $\hat \mu:M\ra \RR^+$
so that $p^*\hat g=p^*(\hat\mu\cdot g)$,
we have $\hat\mu\cdot g=\hat g$.
The conformal class of $g$ contains a l.c.K. metric with parallel Lee 
form.

This finishes the proof of Theorem A.

\vskip0.1cm
Using this theorem,  as is noted in Remark \ref{ex},
we can make clear the case
when the parallel Lee field does not generate a (free) circle action.
 Equivalently,
the parallel Lee field is not a quasi-regular
field.
\begin{cb} [The structure of Vaisman manifolds] \label{base}
Let $(M,g,J)$ be a compact locally conformally K\"ahler, non-K\"ahler manifold 
with
parallel Lee form $\theta$.
\par\noindent{\bf (I)}\ Suppose that the Lee field $\theta^{\#}$
does not generate a circle. Then there exists a l.c.K.
metric $\hat g$ with paralell Lee form in the conformal class of $g$
on $(M,J)$ such that:
\begin{enumerate}
\item The paralell Lee field ${\hat\theta}^{\#}$ of $\hat g$
generates a circle action $S^1$.

\item $(M,\hat g,J)$ is isometric to the quasi-regular
Vaisman manifold\\ $\displaystyle S^1\mathop{{\times}}_{Q}^{} W$ 
where $\displaystyle S^1\ra M\stackrel {\hat\pi}{\lra} W/Q$ is a
Seifert fiber space over a Sasakian orbifold $W/Q$.

\item The K\"ahler form on the universal covering space
$\tilde M$ is identified with $-d(t\pi^*\tilde\eta)$
for which the contact form $\eta$ is a $Q$-invariant
pseudo-Hermitian structure on $W$. Here $t$ is the coordinate for $\RR^+$.
The projection $\pi:\tilde M\ra W$ is the lift of $\hat\pi$ to
$\tilde M=\RR^+\times W$.
\end{enumerate}

\par\noindent{\bf (II)}\ Suppose that the Lee field $\theta^{\#}$
generates a circle $S^1$. Then
$(M,g,J)$ is itself of the quasi-regular
form $\displaystyle S^1\mathop{{\times}}_{Q}^{} W$
for which the K\"ahler form on $\tilde M$ is
$d(t\pi^*\tilde\eta)$ and:
\begin{enumerate}
\item The paralell Lee form $\theta$ lifts to the form
$dt$ on $\tilde M=\RR^+\times W$ where
$t$ is the coordinate for $\RR^+$.
\item The contact form $\tilde\eta$ 
is a $Q$-invariant Sasakian structure on $W$, $i.e.,$  $\pi^*\tilde \eta(X)=
\tilde g(J\xi,X)$ for $X\in\xi^{\#}$.
\end{enumerate}
\end{cb}

\begin{proof}
Let $\mathop{\ISO}(M,g)$ be the group of all
holomorphic isometries of $M$.
As the Lee field $\theta^{\#}$ is Killing, it generates a
$1$-parameter subgroup $\{\varphi_t\}_{t\in\RR}$ of (holomorphic)
isometries on $M$.  Then 
$\{\varphi_t\}_{t\in\RR}\subset\mathop{\ISO}(M,g)$.
Since $\mathop{\ISO}(M,g)$ is compact, the closure $\mathcal T$
of $\{\varphi_t\}_{t\in\RR}$ in $\mathop{\ISO}(M,g)$ is a $k$-torus 
$(k\geq 1)$.\\

{\bf I.}\ When the $1$-parameter group $\{\varphi_t\}_{t\in\RR}$ is 
not a circle,
we can find a sequence of circles $\{S_i^1\}_{i\in \NN}$ which 
approaches
to $\{\varphi_t\}_{t\in\RR}$. Denote by $\{\xi^i\}$ the vector field 
induced by
$S_i^1$ for each $i$. Then there exists $\ell$ with
\begin{equation}\label{xi}
\theta(\xi_\ell)\neq 0\ \mbox{\  everywhere in}\  M.
\end{equation}
To prove this, suppose that  for all $i$ there exists
a sequence of points $\{x_i\}\subset M$  such that 
$\theta_{x_i}(\xi_{x_i})=0$.
As $M$ is compact, there exists a point $x\in W$ such that
$\displaystyle x=\mathop{\lim}_{i\ra \infty}x_i$.
If we note that $\xi_{x_i}$ converges to $(\theta^{\#})_x$, then
$\displaystyle 0=\mathop{\lim}_{i\ra \infty}\theta(\xi_{x_i})
=\theta((\theta^{\#})_x)=g_x(\theta^{\#},\theta^{\#})$, which is impossible.

Now, let $\{\phi_t\}_{t\in\RR}$ (respectively  $\tilde\xi_\ell$)
be the lift of the circle $S^1_{\ell}$
(respectively  $\xi_\ell$ ) to $\tilde M$.
Recall the K\"ahler form
$\Omega=e^{-\tau}\cdot \tilde\omega$ on $\tilde M$ where
$d\tau=\tilde \theta$ (cf. \eqref{po}). 
Note that each $\phi_t$ leaves invariant
$\tilde \omega$ and so it satisfies:
\[
\phi_t^*\Omega=e^{-(\phi_t^*\tau-\tau)}\cdot \Omega.
\]As $\Omega$ is K\"ahler,
$e^{-(\phi_t^*\tau-\tau)}$ is constant on $\tilde M$.
Suppose that every $\phi_t$ is an isometry, $i.e.$
$\phi_t^*\tau-\tau=0$ for all $t$.
Then,
\begin{equation*}
0=\Ll_{\tilde\xi_\ell}(\tau)=d\tau(\tilde\xi_\ell)=
\tilde \theta(\tilde \xi_{\ell})=
\theta(\xi_\ell)\neq 0 \text{\ by}\ \eqref{xi},
\end{equation*}being a contradiction, and so some 
$\phi_t$ is a nontrivial homothety; 
$\phi_t^*\Omega=e^{ct}\cdot \Omega$ for a nonzeo constant $c$.
Hence every $\phi_t$ is a nontrivial homothetic transformation.
By Lemma \ref{ess}, the circle $S^1_\ell$ is purely conformal.
Then the results {\bf 1,2,3} of {\bf I} follow from Theorem A.

We prove {\bf II}. Let $\{\varphi_t\}_{t\in\RR}$ be a lift of $S^1$
to $\tilde M$ which induces a vector field  $\xi$.
Noting that $\varphi_t$ is  an isometry of $\tilde g$ because $\theta$
is parallel,
we can write  
$\varphi_t^*\Omega=e^{-(\varphi_t^*\tau-\tau)}\cdot \Omega$
where $\Omega=e^{-tau}\cdot \tilde \omega$.
Since  $e^{-(\varphi_t^*\tau-\tau)}$ is constant on $\tilde M$ for each $t$,
we may put $\tau\circ \varphi_t-\tau=c\cdot t$ for some constant $c$.
By the definition $\theta(X)=g(\theta^{\#},X)$, so is its lift
$\tilde \theta(\tilde X)=\tilde g(\xi,\tilde X)$ on $\tilde M$.
As $\tilde \theta =d\tau$,
\begin{equation} \label{lift}
0<\tilde g(\xi,\xi)=\tilde\theta(\xi)=\Ll_\xi\tau=
\lim_{t\rightarrow 0}\frac{\tau\circ\varphi_t-\tau}{t}=c.
\end{equation} We may normalize $c=1$ so that
(cf. \eqref{rho}):
\begin{gather} 
\tau\circ \varphi_t(w)-\tau (w)=t \text{\ for }\ w\in W,\label{con}\\
\varphi_t^*\Omega=e^{-t}\cdot \Omega \text{\ on} \tilde M.\label{rrho}
\end{gather}
In particular, the group
$\{\varphi_t\}_{t\in\RR}$ is isomorphic to $\RR$.
As in the proof of Theorem A, we have a contact form $\eta$ on $W$.
However, comparing \eqref{rrho} with \eqref{rho}, this will be
defined as
\begin{equation}\label{ttrei}
\eta_w(X_w)=e^{t}\cdot\Omega(\tilde X_{\f_tw},\xi_{\f_tw})
\mbox{\ i.e.,}\ \pi^*\eta=-e^t\cdot \iota_\xi\Omega,
\end{equation} where
 $\tilde X_{\f_tw}\in T_{\f_tw}\tilde M$ such that $\pi_*\tilde
X=X$. (Compare \eqref{trei}.)
Noting also that $\Ll_\xi\Omega=-\Omega$ in this case,
we have that (cf. \eqref{eta})
\begin{equation}\label{eeta}
d(e^{-t}\cdot\pi^*\eta)=\Omega \;\, \text{on}\  \tilde M.
\end{equation}

On the other hand, we can define a $1$-form $\tilde\eta$ on $W$
({\em that is, a Sasakian contact structure}, see \cite{DO}) by
\[
\tilde\eta(X)=\tilde g(J\xi,\tilde X)
\] for $\tilde X\in\xi^{\perp}$ with respect to $\tilde g$
such that $\pi_*\tilde X=X$. If we note
$h=e^{-\tau}\cdot \tilde g$, then the distribution $\xi^{\perp}$
is the same as that of $h$. As $\pi_*:
\xi^{\perp}\ra TW$ is isomorphic, $\tilde \eta$ is well defined.\\

We compare these contact forms first. By $\Omega=e^{-t}\cdot \tilde\omega$
and \eqref{ttrei}, we compute
\[
\pi^*\tilde\eta(\tilde X)=-\tilde g(\xi,\tilde X)
=-\tilde\omega(\xi,\tilde X)=-e^{\tau}\cdot\Omega(\xi,\tilde X)
=e^{\tau}\cdot(e^{-t}\cdot \pi^*\eta)(\tilde X)
\] so that
\begin{equation}\label{kan}
\pi^*\tilde\eta=e^{\tau-t}\cdot \pi^*\eta \text{\ on}\ \tilde M.
\end{equation}
Moreover, $\varphi_s^*(\tau-t)=\tau\circ \varphi_s-(t+s)=
\tau-t$ on $W$ from \eqref{con}. Thus there is a map $\nu
: W \ra \RR$ with $\pi^*\circ \nu=\tau-t$. Hence,
\eqref{kan} implies:
\begin{equation}\label{kankan}
\tilde\eta=e^{\nu}\cdot \eta \text{\ on}\ W.
\end{equation}

Finally, recall from Lemma 
\ref{s1} that
$1=s(w)=\Omega(J\xi_w,\xi_w)=e^{-\tau(w)}\cdot
\tilde \omega(J\xi_w,\xi_w)=
e^{-\tau(w)}\cdot\tilde g(\xi_w,\xi_w)=e^{-\tau(w)}$
because $c=1$ in \eqref{lift}.
Thus, from   equation \eqref{con}
we derive:
\begin{equation} \label{st}
e^{\tau\circ \varphi_t(w)}=e^t \text{\ for all}\ \varphi_t(w)\in \tilde M.
\end{equation}
Using \eqref{st} and \eqref{eeta}, \eqref{kan}
\begin{eqnarray*}
\tilde\omega&=&e^{\tau}\cdot\Omega=e^{\tau}\cdot
d(e^{-t}\pi^*\eta)\\
&=&e^{\tau}\cdot d(e^{-t}e^{t-\tau}\cdot \pi^*\tilde\eta)=
e^{\tau}\cdot d(e^{-\tau}\cdot \pi^*\tilde\eta)\\
&=&e^t\cdot d(e^{-t}\cdot \pi^*\tilde\eta)
\end{eqnarray*} on $\tilde M=\RR^+\times W$.
Therefore $\tilde g$ has the parallel Lee form
$dt$ such that
$p^*\theta=dt$
for $\theta$ on $M$.
This proves {\bf II.}

\end{proof}

\begin{re}
\ {\bf (1)}\
When a compact l.c.K. manifold $(M,g,J)$ has
the quasi-regular form $S^1\times W/Q$ in the conformal class
of $g$,  it is still vague whether or not $g$ itself has parallel Lee
form.\\
\smallskip 
\noindent
{\bf (2)}\ The Inoue
surface $S_N^-$ has a l.c.K. metric without parallel Lee form
(see \cite{Tr}). The associated Lee field $\theta^{\#}$
is Killing (without constant norm) and generates an $S^1$-action by
holomorphic isometries. It acts freely so that the orbit space is the
$3$-dimensional solvmanifold. On the other hand, a compact
Sasakian manifold admits a nontrivial $T^k$-action $(k\geq 1)$
generated by the Reeb field. If the Sasakian manifold
is a closed aspherical manifold, then the fundamental group
has a nontrivial center (at least, containing
a free abelian group of rank $k$.)
As the $3$-dimensional solvmanifold has no center, it admits
no Sasakian structure.
Hence, this $S^1$ is not purely conformal. Note that
the one-parameter subgroup of holomorphic transformations
generated by $J\theta^{\#}$ is not a circle
and it does not leave invariant the fundamental two-form.
\end{re}
\section{Lee-Cauchy-Riemann transformations}
In this section, we  consider another kind of transformations of
 l.c.K. manifolds.
As above,
there exists an orthonormal,
local coframe\\
$\{\theta, \theta\circ J, \theta^\al,\ov{\theta}^\al\}_{\al=1,
\cdots, n-1}$ adapted to a l.c.K. manifold $(M,g,J)$. Here we have
 =
$\displaystyle\omega|\{\theta^{\#},J\theta^{\#}\}^{\perp}=\mathop{\sum}_{1}^{n-1}
\delta_{\al\be}\theta^\al\wedge\ov{\theta}^\be$.
 Consider the diffeomorphism $f$ which transforms:
\begin{equation} \label{rel}
\begin{split}
f^*\theta=\theta,\ \  f^*(\theta\circ J)=\lambda\cdot(\theta\circ 
J),& \\
f^*\theta^{\al}=\sqrt \lambda\cdot\theta^\be U^{\al}_{\be}+
(\theta\circ J)\cdot v^{\al},& \\
f^*{\bar \theta}^{\al}=\sqrt \lambda\cdot
{\bar\theta}^\be\ov{U}^\al_\be+(\theta\circ J)\cdot\ov{v}^{\al} &
\end{split}
\end{equation}
for some positive, smooth function $\lambda$, and
a matrix $U^\al_\be\in {\rm U}(n-1)$.

Let  $\{\theta^{\#}, J\theta^{\#}\}^{\perp}$
be the subbundle whose vectors are perpendicular
to $\{\theta^{\#}, J\theta^{\#}\}$ with respect to $g$.
As $\{\theta^{\#}, J\theta^{\#}\}^{\perp}$ coincides with the
distribution
 $\mathop{\Null}\, \theta\cap\mathop{\Null}\, \theta\circ J$,
it is invariant under $J$.
By the above definition,
$f$ maps the subbundle $\{\theta^{\#}, J\theta^{\#}\}^{\perp}$onto 
itself
and is holomorphic on it, $i.e.,$ $f_*\circ J=J\circ f_*$.

On the other hand, since $f^*\theta=\theta, f^*(\theta\circ 
J)=\lambda\cdot
(\theta\circ J)$,
$f_*\theta^{\#}=\theta^{\#}+A,
 f_*J\theta^{\#}=\lambda\cdot J\theta^{\#} +B$
for some $A,B\in\{\theta^{\#}, J\theta^{\#}\}^{\perp}$.
Hence $TM$ decomposes into the direct sum 
$\{f_*\theta^{\#},f_*J\theta^{\#}\}
\oplus\{\theta^{\#}, J\theta^{\#}\}^{\perp}$.
Define an almost complex structure $J'$ on $M$ simply
to be
\begin{eqnarray*}
J'f_*\theta^{\#}&=&f_*(J\theta^{\#}),\\
J'f_*(J\theta^{\#})&=&-f_*\theta^{\#},\\
J'|\{\theta^{\#}, J\theta^{\#}\}^{\perp}&=&J.
\end{eqnarray*}
If a vector $X=a\theta^{\#}+b(J\theta^{\#})+V\in TM$
for $V\in\{\theta^{\#}, J\theta^{\#}\}^{\perp}$, then
we can check that $f_*(JX)=J'f_*(X)$. Thus $f$ is $J'$-holomorphic on 
$M$
so that $J'$ is a complex structure.\\

We put
$\omega'={f^{-1}}^*\omega$. Since ${f^{-1}}_*\circ J'=J\circ 
{f^{-1}}_*$,
it is easy to see that the two-form
$\omega'$ is $J'$-invariant.
Letting $\theta'={f^{-1}}^*\theta$,
 $d\omega'=\theta'\wedge \omega'$ such that $d\theta'=0$.
If we set $g'(X,Y)=\omega'(J'X,Y)$, then
$g'$ is $J'$-invariant l.c.K. metric on $M$.
We can prove that $g'(f_*X,f_*Y)=g(X,Y)$.
Hence $f$ defined by \eqref{rel}
is a l.c.K. diffeomorphism of $(M,g,J)$ onto
$(M,g',J')$. We call such a transformation $f$ a
\emph{Lee-Cauchy-Riemann $($LCR$)$
transformation}  from one l.c.K. manifold $(M,g,J)$
to another l.c.K. manifold $(M,g',J')$.

If $f'$ is a $LCR$ transformation satisfying the equations that
${f'}^*\theta=\theta$, ${f'}^*(\theta\circ J)=\lambda'\cdot(\theta\circ J)$,
${f'}^*\theta^{\al}=\sqrt {\lambda'}\cdot\theta^\be {U'}^{\al}_{\be}+
(\theta\circ J)\cdot {v'}^{\al}$,\\ ${f'}^*{\bar \theta}^{\al}=\sqrt {\lambda'}
\cdot{\bar\theta}^\be_\al\ov{U'}^\al_\be+(\theta\circ J)\cdot\ov{v'}^{\al}$,
then it is easy to see that
\begin{equation} \label{rrel}
\begin{split}
& {(f'\circ f)}^*\theta=\theta, \ \ {(f'\circ f)}^*(\theta\circ J)=
({f^*\lambda'}\cdot \lambda)\cdot(\theta\circ J),\\
&{(f'\circ f)}^*\theta^{\al}=
\sqrt {f^*\lambda'\cdot\lambda}\cdot\theta^\be W^{\al}_{\be}+(\theta\circ J)
\cdot w^{\al},\\
& {(f'\circ f)}^*{\bar \theta}^{\al}=
\sqrt {f^*\lambda'\cdot\lambda}\cdot{\bar\theta}^\be\ov{W}^\al_\be+
(\theta\circ J)\cdot\ov{w}^{\al}.
\end{split}
\end{equation}The composition
$f'\circ f$ is also a $LCR$ transformation from $(M,g,J)$
to another l.c.K. manifold.
Denote by $\mathop{\Aut}_{LCR}(M,g,J,\theta)$
 (or simply by $\mathop{\Aut}_{LCR}(M)$)
the group of all Lee-Cauchy-Riemann transformations
on a l.c.K. manifold $(M,g,J)$ adapted to the Lee form $\theta$.

Consider the subgroup $G$ of $\mathop{\GL}(2n,\RR)$
consisting of the following elements:
$$
\left\{\left(\begin{array}{cccccc}
1&0&\ & 0        &\  &0\\
0&u&\ &v^{\alpha}&\ &v^{\bar\alpha}\\
0&0&\ &\sqrt u U_{\beta}^{\alpha}&\ &0\\
0&0&\ &0&\ &\sqrt u U_{\bar\beta}^{\bar\alpha}
\end{array}\right)
\ |\ u\in\RR^+,v^{\alpha}\in\CC, U_{\beta}^{\alpha}\in {\rm 
U}(m)\right\},
$$where $m=n-1$. Let $G\ra P\ra M$  be the principal bundle of
the $G$-structure consisting of the above coframes
$\{\theta, \theta\circ J, \theta^{\alpha},\theta^{\bar\alpha}\}$.
If we note that $G$ is isomorphic to the semidirect
product $\CC^n\rtimes ({\rm U}(m)\times \RR^+)$,
then the Lie algebra ${\mathfrak g}$
is isomorphic to $\CC^n+{\mathfrak u}(m)+\RR$.
In particular, the matrix group
${\mathfrak g}\subset {\mathfrak gl}(2n,\RR)$
has no element of ${\rm rank}$
 1, \emph{i.e.} it is  \emph{elliptic} (note that
$\CC^n$ is of infinite type, while ${\mathfrak u}(m)+\RR$
is of order $2$ (cf. \cite{Ko}).) As $M$ is assumed to be
compact, the group of automorphisms $\mathcal U$ of
$P$ is a (finite dimensional) Lie group.
Since $\mathop{\Aut}_{LCR}(M)$
is a closed subgroup of $\mathcal U$,
$\mathop{\Aut}_{LCR}(M)$ is a Lie group.
The compactness of the group $Aut_{LCR}(M)$,
on a compact l.c.K. manifold $(M,g,J)$
will be discussed in a subsequent paper.
By Proposition \ref{lck},
 each element of $Aut_{l.c.K.}(M)$
satisfies
$f^*\theta=\theta$, $f^*(\theta\circ J)=(\theta\circ J)$
and $f^*\omega=\omega$ with respect to some specific coframes
$\{\theta, \theta\circ J, \theta^{\alpha},\theta^{\bar\alpha}\}$.
The fact that  $f^*\omega=\omega$ implies that
$f^*\theta^{\al}=\theta^\be U^{\al}_{\be}$,
$f^*{\bar \theta}^{\al}={\bar\theta}^\be\ov{U}^\al_\be$
for some positive matrix $U^\al_\be\in {\rm U}(n-1)$.
Note that the elements of $Aut_{l.c.K.}(M)$ are
also viewed as $LCR$ transformations.

When a noncompact Lee-Cauchy-Riemann
transformation subgroup of $Aut_{LCR}(M)$ acts on a
compact l.c.K. mamifold $M$,
we prove a similar property
to the noncompact $CR$-action
on a compact $CR$-manifold.
The next result,
the proof of which will be mostly situated in the realm of CR-geometry,
characterizes the Hopf manifolds, up to biholomorphism, among the
compact locally conformally K\"ahler manifolds.

\begin{tc}
Let $(M,g,J)$ be a compact locally conformally K\"ahler manifold with
Lee form $\theta$. 
Suppose that there exists a
closed noncompact subgroup $\RR=\{\h_t\}_{t\in\RR}$
of $Aut_{LCR}(M)$ satisfying that
\begin{enumerate}
\item $(\h_{t})_*(\theta^{\#})=\theta^{\#}$.
\item $(\h_t)_*$ preserves the distribution
$\{\theta^\sharp, (\theta\circ J)^\sharp\}^\perp$
and is holomorphic on it.
\end{enumerate}
Then the following hold:
\par\noindent {\bf (i)}\ If
$\theta$ is parallel and the Lee field $\theta^{\#}$ generates a
circle $S^1$, then
$(M,g,J)$ is holomorphically isometric
 to a finite quotient of a Hopf manifold
$S^1\times S^{2n-1}/\Gamma$ where $\Gamma\subset
S^1\times{\rm U}(n-1)$.
\smallskip
\par\noindent{\bf (ii)}\ If the Lee field $\theta^{\#}$
generates a purely conformal circle $S^1$, then
 $(M,g,J)$ is holomorphically conformal to 
a finite quotient $S^1\times S^{2n-1}/\Gamma$ $(\Gamma\subset
S^1\times{\rm U}(n-1))$.
\end{tc}
From this theorem, we obtain the following:
\begin{cd}
Let $(M,g,J)$ be a compact locally conformally K\"ahler manifold which
 admits a closed subgroup $\CC^*$
of Lee-Cauchy-Riemann transformations adapted
to the Lee form $\theta$. Then,
\par\noindent{\bf (1)}\ If $\theta$ is parallel and a subgroup
$S^1$ of $\CC^*$ induces the Lee field $\theta^{\#}$, 
then some finite cover of $M$ is holomorphically isometric with
a Hopf manifold $S^1\times S^{2n-1}$.
\par\noindent
In general,
\par\noindent{\bf (2)}\
If there exists a purely
conformal circle $S^1$ in $\CC^*$ which induces the Lee field
$\theta^{\#}$, then
$M$ is holomorphically conformal to
a finite quotient of a Hopf manifold $S^1\times S^{2n-1}$.
\end{cd}

\subsection*{Proof of theorem C}\hfill 
\par Let $(M,g,J)$ be a compact locally conformally K\"ahler manifold
satisfying the conditions {\bf 1,\ 2}.
We work with the same l.c.K. metric $g$ on $(M,J)$ for the case {\bf (i)},
while applying Theorem A for the case {\bf (ii)} first of all,
there is a l.c.K. metric $\hat g$
with parallel Lee form which is conformal to $g$. (As a consequence,
note that the Lee field $\theta^{\#}$ of $g$ becomes  Killing
with respect to $\hat g$.)
Then we work with $(M,\hat g,J)$ to show that it is
holomorphically isometric with a finite quotient of a Hopf manifold.
(As a matter of fact, $(M,g,J)$ is holomorphically conformal
to a finite quotient of a Hopf manifold.).

We retain the same notation $\h_{t}$  for the lift of $\RR=<\h_t>$
to the universal covering $\tilde M$. Suppose that the hypothesis
of {\bf (i)} is satisfied.
Let $\hat\varphi_t$ be an element of
$S^1$ generated by the Lee field $\theta^{\#}$. By {\bf 1}, each $\h_t$ commutes with every element $\hat\varphi_t$.
Then each lift $\h_t$ to $\tilde M$ also commutes with the elements of
the lift $\RR^+$ of $S^1$. 
Let $\RR^+=\{\varphi_t\}_{t\in\RR}$ be a $1$-parameter group.
(See {\bf (II)} of Corollary B.)
Then
$$
\h_{t}\circ \f_s=\f_s\circ\h_t\ \ (t,s\in\RR).
\leqno{(*)}
$$
On the other hand, if the hypothesis
of {\bf (ii)} is satisfied, then  as in the proof of Theorem A, the
$S^1$-action lifts to
a $\RR^+ (=\{\varphi_t\}_{t\in\RR})$-action which satisfies 
$(*)$.

In each case, let $\xi$ be a vector field on $\tilde M$ induced by
the $\RR^+$-action and note that $p_*\xi=\theta^{\#}$. By $(*)$,
the projection $\displaystyle\RR^+\ra \tilde M\stackrel{\pi}{\lra}W$
induces a $\RR=\{\h_t\}$-action on $W$.
We have the following commutative diagram:
\begin{equation}\label{com}
\begin{CD}
\mathbb{Z}@>>>\pi_1(M)@>>>Q\\
@VVV          @VVV        @VVV\\
\RR^+@>>>(\RR^+\times \RR, \tilde M)@>\pi >>(\RR, W)\\
@VVV @VVpV @VVpV\\
S^1 @>>> (S^1\times \RR, M) @>\pi^* >> (\RR, M^*)
\end{CD}
\end{equation}
where $M^*=M/S^1$ and in the top line of the diagram,
$Q=\pi_1(M)/\mathbb{Z}$.
By {\bf 2}, each $\h_t$ preserves the distribution
$\{\theta^\sharp, (\theta\circ J)^\sharp\}^\perp$.
We observe that $p_*\xi=\theta^{\#}$,
$p_*(J\xi)=J\theta^{\#}=-(\theta\circ J)^{\#}$, hence  
$p$ maps $\{\xi,J\xi\}^{\perp}$ isomorphically onto 
$\{\theta^\sharp, (\theta\circ J)^\sharp\}^\perp$ at each point of $M$.
Since $\pi_* :\{\xi,J\xi\}^{\perp}
\ra \mathop{\Null}\, \eta$ is $J$-isomorphic,
the induced $\RR$-action $<\h_t>$ on $W$
preserves $\mathop{\Null}\, \eta$
on which it is holomorphic,
$i.e.,$ $\h_{t}^*\eta=\lambda_t\cdot \eta$ for some
map $\lambda_t:W\ra \RR^+$, and ${\h_t}_*\circ J=J\circ{\h_t}_*$.
Hence $\RR=\{\h_t\}$ is a closed subgroup of $CR$-transformations
of $(W,\eta,J)$.
Note that the contact structure $\eta$ of the Sasakian manifold $W$ 
gives a strictly pseudoconvex, pseudo-Hermitian structure 
and the quotient space $W/Q$ inherits a $CR$-structure
from $(\eta,J)$ by the property {\bf 1} of \eqref{inv}. Then
the $CR$-action of $W$ induces a $CR$-action of $\RR$ on $W/Q$.

When $W/Q$ happens to be a compact smooth manifold, we can apply 
Webster's theorem in \cite{W} to yield that $W/Q$ is spherical.
But in general, $W/Q$ is a compact orbifold. Even so,
we can still show that it is spherical.
The key step will be, as in \cite{W}, to prove:
\begin{pr}\label{sf} The $CR$-manifold
$(W,\{\mathop{\Null}\, \eta,J\})$ is spherical $($or CR-flat, in 
another terminology,
 meaning that its Chern-Moser curvature tensor $S$
vanishes identically on $W)$.
\end{pr}
\begin{proof}
Suppose $W$ is not spherical. Then the set $V=\{x\in W\ |\ S_x\neq 
0\}$ is a
non-empty, open subset of $W$. As $S$ is a $CR$ invariant (see \emph{e.g.}
 \cite{W}), the set $V$ is preserved by any CR-automorphism of $W$.
If a contact form $\eta$ is replaced by $\eta'=u\cdot \eta$, then the 
norm
of the Chern-Moser tensor with respect to the Levi form associated to 
$\eta$
 satisfies the equality $\Vert S\Vert_\eta=u\cdot \Vert S\Vert_{\eta'}$ (cf. 
\cite{W1}). Choose $u=\Vert S\Vert_\eta$ on $V$. Then we obtain a
pseudo-Hermitian structure $(\eta',J)$ defined on $V$
such that
\begin{equation}\label{nor}
\Vert S\Vert_{\eta'}=1\ \text{on $V$}.
\end{equation}
If $f\in{Aut}_{CR}(W)$, then $f^*\eta'=\lambda\cdot \eta'$.
The above relation shows that
$\Vert S\Vert_{\eta'}(x)=\lambda\cdot\Vert S\Vert_{f^*\eta'}(x)
=\lambda\cdot\Vert S\Vert_{\eta'}(f(x))$.
By \eqref{nor}, $\lambda=1$ on $V$.
Then we can check that
for each element $f$ of $Aut_{CR}(W)$, the transformation \eqref{rel}
can be reduced to the following:
\begin{eqnarray*}
f^*\eta' &=& \eta' \\
f^*{\theta'}^\al &=&{\theta'}^\be U^\al_\be\\
f^*{\overline{\theta}'}^{\al} &=& {\overline{\theta}'}^\be
{\ov{U}}^\al_\be\\
d\eta'&=&\mathop{\sum}_{\al,\be}^{}{\theta'}^{\alpha}
\wedge{\overline{\theta}'}^\be.
\end{eqnarray*}Thus the bundle of such coframes $\{\eta', \eta'\circ J,
 {\theta'}^\al,\ov{\theta'}^\al\}_{\al=1,\cdots, n-1}$  with
$d\eta'=\mathop{\sum}_{\al,\be}^{}{\theta'}^\al\wedge
{\overline{\theta}'}^\be$ gives rise to a principal bundle,
restricted to $V$:
\begin{equation}\label{bun}
{\rm U}(n)\ra P'\stackrel{q}{\lra} V,
\end{equation}
for which we note that
\begin{equation}\label{closed}
\text{$Aut_{CR}(W)$ is a \emph{closed} subgroup of $Aut_{U(n)}(V)$.}
\end{equation}
On the other hand, as ${\rm U}(n)$ is of order $1$ (as a subgroup of
 ${\rm O}(2n)$),
 the manifold $P'$ has a $\{1\}$-structure. Hence, by Theorem 3.2 in 
\cite{Ko},
 for any fixed $u'\in P'$, the orbit map $Aut_{U(n)}(P')\ra P'$
is a proper embedding, and the orbit
 $Aut_{U(n)}(P')\cdot u'$ is closed in $P'$. As any automorphism of $W$
produces, by differentiation, an automorphism of $P'$,
we may consider the subgroup
$Aut_{U(n)}(V)_*=\{df\ |\ f\in Aut_{U(n)}(V)\}$. Restricting to
CR-automorphisms of $W$ and using \eqref{closed}, we derive that the 
orbit
 $Aut_{CR}(W)_*\cdot u'$ is closed in $P'$ and that
$Aut_{CR}(W)_*\cdot u'$ is homeomorphic to
$Aut_{CR}(W)$.  Now recall that we have a closed
subgroup $\RR$ in $Aut_{CR}(W)$. Denote with $\RR_*$ its lift to $P'$
 ($\RR_*$ contains the differentials of the
 flows associated to the $\RR$-action). Hence, $\RR_*\cdot u'$ is closed 
in
 $P'$. Now, if we can prove that this orbit is contained in a certain 
compact
 subset $K'\subset P'$, we see that $\RR_*\cdot u'$ has to be compact, 
in
 contradiction with it being homeomorphic with $\RR$ by the proper 
embedding.

To construct $K'$, set $x'=q(u')\in V$ so that
$q(\RR_*\cdot u')=\RR\cdot x'$.
From the diagram \eqref{com}, $p$ is a quotient (continuous) map and $W$ 
is
locally compact (as a manifold), $M^*$ is compact, hence we may find a
connected, compact subset $C\subset W$ such that $x'\in C'$ and 
$p(C)=M^*$.
Hence $\displaystyle W=\mathop{\cup}_{\alpha\in Q}\alpha\cdot C$.
In particular, as $\RR\cdot x'\subset V$, we obtain
$\displaystyle\RR\cdot x'\subset \mathop{\cup}_{\alpha\in Q}
\alpha\cdot C$. But $Q$ acts properly discontinuously,
so only a finite number of
translated $\alpha\cdot C$ meet $C$. Since $\RR\cdot x'$ is connected, 
we may
write $\displaystyle\RR\cdot x'\subset 
C':=\mathop{\cup}_{i=1}^k\alpha_i
\cdot C$ and note that $C'$ is
compact. Hence, the closure $\ov{\RR\cdot x'}$ is compact.
But \emph{a priori}, it might exit $V$.
We prove that this is not the case and, in fact,
$\ov{\RR\cdot x'}\subset V$ following \cite{W}.
Then the inverse image $K'=q^{-1}(\ov{\RR\cdot x'})$
of the bundle \eqref{bun} is the one desired.

To this end,
let $A$ be the vector field on $W$ induced by the considered 
$\RR$-action
$<\h_t>$.
We first prove:
\begin{equation}\label{A}
 \eta'(A) \; \text{does not vanish identically on $V$.}
\end{equation}
Indeed, by absurd, if $\eta'(A)=0$ on $V$, then $A\in \mathop{\Null}\, 
\eta'$ on $V$.
For any $X\in\mathop{\Null}\, \eta'$, as $A$ an infinitesimal
contact transformation, we have $0=(\Ll_A\eta')(JX)=A(\eta'(JX))-\eta'([A,JX])$ = $2d\eta'(A,JX)$,
in contradiction with $(V, \eta')$ being strictly pseudo-convex,
its Levi form
$d\eta'$ is positive definite, in particular non-degenerate.
If we note that $x'\ (=q(u'))$ can be chosen arbitrary in $V$, we may 
suppose
$\eta'(A_{x'})=d$ for some suitable, fixed $d\in \RR^*$.

We consider the non-empty set
$D=\{x\in V\; ;\; \eta'(A_x)=d\}$ and show that:
\begin{equation}\label{D}
\text{$D$ is closed in $W$.}
\end{equation}
Here is a simple argument of general topology.
Observe first that, by  $\eta'=u\cdot\eta$ with $u=\Vert 
S\Vert_\eta$, $u(x)$
tends to $0$ when $x$ approaches the boundary of $V$ in $W$. Now let 
$x\in
\overline{D}$ and choose a sequence $\{x_n\}\subset D$
which converges to $x$.
We have $d=\eta'_{x_n}(A_{x_n})=u(x_n)\cdot\eta_{x_n}(A_{x_n})$. But
$\eta_{x_n}(A_{x_n})\rightarrow \eta_x(A_x)<\infty$ as $x_n\rightarrow 
x$,
 hence $u(x_n)\not\rightarrow 0$. This means that $x$ is not a boundary 
point
of $V$, thus $x\in V$. This implies $\eta'(A_x)=\lim 
\eta'(A_{x_n})=d$,
yielding $x\in D$ and proving \eqref{D}.

Since $f^*\eta'=\eta'$ for $f\in Aut_{CR}(W)$,
the equality $\eta'(\h_{t*}A_{x'})$ = $\eta'(A_{x'})$ = $d$ implies
$\RR\cdot x'\subset D$. Obviously,
$\ov{\RR\cdot x'}\subset V$.
This ends the proof of Proposition \ref{sf}.
\end{proof}
\par Now, in order to finish the proof of Theorem C, we analyse the following 
diagram:
\begin{equation*}
\begin{CD}
\, Q\\
@VVV\\
(\RR, W)@>(\rho,\, dev)>>(P\mbox{\rm U}(n,1), S^{2n-1})\\
@VVV\\
(\RR, W/Q)
\end{CD}
\end{equation*}
where in the horizontal line we have the developing pair (see \cite{Ku} 
for the
definition) of $(\RR, W)$ and  $P\mbox{\rm U}(n,1)$ is the group of 
CR-automorphisms of
the standard sphere viewed as boundary of the complex hyperbolic space
(cf. \cite{Kam}). If $W$ is not CR-equivalent with $S^{2n-1}$, then
$X:=dev(W)$ may consist of $S^{2n-1}-\{\infty\}$ or 
$S^{2n-1}-\{0,\infty\}$
 (see Theorem 4.4 in \emph{loc. cit.}). We show that these two cases 
lead to
contradiction. Indeed, if by absurd either of the two cases occur,
we may pull back by $dev$ their metrics to $W$ obtaining a $Q$-invariant
metric. This one descends to a complete metric of the compact orbifold 
$W/Q$.
We take the associate distance function (which is continuous)
and we can lift
it back to $W$ to a complete distance
(all this was necessary because $W$ is not compact!).
It follows that $dev$ is a covering map, hence a diffeomorphism
 on the image. As a consequence, a finitely generated
subgroup $Q$ is isomorphic with $\rho(Q)$ in $P\mbox{\rm U}(n,1)$.
 According to Selberg's result (cf. \cite{Ra}, Corollary 6.14)
there exists a subgroup $Q'$ of finite index in $Q$ which is torsion 
free.
Hence, $W/Q'$ is a compact manifold.  Then the group of
$CR$ automorphisms of $W/Q'$ is compact
because it is a closed subgroup of
the group of isometries $W/Q'$ induced by the
metric from $W$. Still, we know that it must contain a closed subgroup 
$\RR$.
This contradiction shows that $dev$ maps $W$ onto $S^{2n-1}$
$CR$-diffeomorphically.
This finishes the proof.
\begin{re}
The $\RR$-action on $S^{2n-1}$ is characterized as
either loxodromic $(=\RR^+)$
or parabolic $(=\mathcal R)$ for which
$\RR^+$ has exactly two fixed points $\{0,\infty\}$ or
$\mathcal R$ has the unique fixed point $\{\infty\}$ on $S^{2n-1}$.
Moreover,
since $\RR$ centralizes $Q$, it implies
either $\RR\times Q\subset \mathcal R\times \mbox{\rm U}(n)$
or $\RR\times Q\subset \RR^+\times \mbox{\rm U}(n)$ where
$Q$ is a finite subgroup by properness. As $W/Q(\approx S^{2n-1}/Q)$
is an orbifold, such a finite subgroup may exist, contrary to the case 
that
$W/Q$ is a compact smooth $CR$-{\it manifold}.
\end{re}

\noindent{\bf Acknowledgement.} The second named author thanks the Department of 
Mathematics of the Tokyo Metropolitan University for support and 
hospitality during July 2000 when part of this paper was written.

\end{document}